\documentclass[12pt]{article}
\usepackage{amssymb}
\usepackage{verbatim}
\usepackage{array}
\usepackage{latexsym}
\usepackage{enumerate}
\usepackage{amsmath}
\usepackage{amsfonts}
\usepackage{amsthm}
\usepackage[english]{babel}

\title{\bf{A note on cyclic semiregular subgroups of some
$2$-transitive permutation groups}}
\date{}

\author{M.~Giulietti  and G.~Korchm\'aros}

\newtheorem{theorem}{Theorem}[section]
\newtheorem{proposition}[theorem]{Proposition}

\newtheorem{claim}[theorem]{Claim}

{\theoremstyle{definition}
\newtheorem*{definition*}{Definition}

\newtheorem*{proposition*}{Proposition}
\newtheorem*{corollary*}{Corollary}
\newtheorem*{lemma*}{Lemma}

\def\cQ{\mathcal Q}
\def\cO{\mathcal O}

\def\cU{\mathcal U}

\def\K{\mathbb{K}}
\def\PG{{\rm{PG}}}

\def\Fq{{\mathbb F}_q}

\def\Fn{{\mathbb F}_n}
\def\fns{{\mathbb F}_{n^2}}

\newcommand{\PSL}{\mbox{\rm PSL}}
\newcommand{\PGL}{\mbox{\rm PGL}}

\newcommand{\PSU}{\mbox{\rm PSU}}
\newcommand{\PGU}{\mbox{\rm PGU}}

\newcommand{\Sz}{\mbox{\rm Sz}}
\newcommand{\Ree}{\mbox{\rm Ree}}

\newcommand{\diag}{\mbox{\rm diag}}

\newcommand{\vf}{\varphi}

\newcommand{\gO}{\Omega}

\newcommand{\ha}{{\textstyle\frac{1}{2}}}

\newcommand{\bA}{{\bf A}}

\newcommand{\bS}{{\bf S}}

\sloppy

\begin{document}
\maketitle

    \begin{abstract}
We determine the semi-regular subgroups of the $2$-transitive
permutation groups $\PGL(2,n),\PSL(2,n), \PGU(3,n),
\PSU(3,n),\Sz(n)$ and $\Ree(n)$ with $n$ a suitable power of a
prime number $p$.
    \end{abstract}\thanks{2000 {\em Math. Subj. Class.}: 14H37 }

\thanks{{\em Keywords}: 2-transitive permutation groups}

    \section{Introduction}
The finite $2$-transitive groups play an important role in several
investigations in combinatorics, finite geometry, and algebraic
geometry over a finite field. With this motivation, the present
notes are aimed at providing some useful results on semi-regular
subgroups of the $2$-transitive permutation groups $\PGL(2,n),\,
\PSL(2,n),\,\PGU(3,n),\,\PSU(3,n),\,\Sz(n)$ and $\Ree(n)$ where
$n$ is a suitable power of a prime number $p$.

\section{The projective linear group}
The projective linear group $\PGL(2,n)$ consists of all linear
fractional mappings,
 $$
\varphi_{(a,b,c,d)}:\quad x\mapsto \frac{ax+b}{cx+d},\quad
ad-bc\neq 0,
$$ with $a,b,c,d\in\Fn$. The order of $\PGL(2,n)$ is
$n(n-1)(n+1)$.

Let $\square$ be the set of all non-zero square elements in $\Fn$.
The special projective linear group $\PSL(2,n)$ is the subgroup of
$\PGL(2,n)$ consisting of all linear fractional mapping
$\varphi_{(a,b,c,d)}$ for which $ad-bc\in \square$. For even $n$,
$\PSL(2,n)=\PGL(2,n)$. For odd $n$, $\PSL(2,n)$ is a subgroup of
$\PGL(2,n)$ of index $2$.

For $n\geq 4$, $\PSL(2,n)$ is a non-abelian simple group. For
smaller values of $n$, $\PGL(2,2)\cong \PSL(2,3)\cong\
{\rm{Sym}}_3$. 
For this reason, we only consider the case of $n\geq 4$.

The above fractional mapping $\varphi_{(a,b,c,d)}$ defines a
permutation on the set $\Omega=\Fn\cup\{\infty\}$ of size $n+1$.
So, $\PGL(2,n)$ can be viewed as a permutation group on $\Omega$.
Such a permutation group is sharply $3$-transitive on $\Omega$, in
particular $2$-transitive on $\Omega$, and it is defined to be the
{\em natural $2$-transitive permutation representation of
$\PGL(2,n)$}. In this context, $\PSL(2,n)$ with $n$ odd can be
viewed as permutation group on $\Omega$. Such a permutation group
is $2$-transitive on $\Omega$, and it is defined to be the {\em
natural $2$-transitive permutation representation of $\PGL(2,n)$}.

The subgroups of $\PSL(2,n)$ were determined by Dickson, see
\cite[Hauptsatz 8.27]{huppertI1967}.
\begin{theorem}
{\em Dickson's classification of subgroups of $\PSL(2,n)$}: If $U$
is a subgroup of $\PSL(2,n)$ with $n=p^r$, then $U$ is one of the
following groups:
\begin{itemize}
\item[{\rm{(1)}}] An elementary abelian $p$-group of order $p^m$
with $m\le r$.

\item[{\rm{(2)}}] A cyclic group of order $z$ where $z$ is a
divisor of $2^r-1$ or $2^r+1$, if $p=2$, and a divisor of $\ha
(p^r-1)$ or $\ha (p^r+1)$, if $p>2$.

\item[{\rm{(3)}}] A dihedral group of order $2z$ where $z$ is as
in {\rm (2)}.

\item[{\rm{(4)}}] A semidirect product of an elementary abelian
$p$-group of order $p^m$ and a cyclic group of  order $t$ where
$t$ is a divisor of $p^{{\rm{gcd}}(m,r)}-1$.

\item[{\rm{(5)}}] A group isomorphic to $A_4$. In this case, $r$
is even, if $p=2$.

\item[{\rm{(6)}}] A group isomorphic to $S_4$. In this case,
$p^{2^r}-1\equiv 0 \pmod {16}$.

\item[{\rm{(7)}}] A group isomorphic to $A_5$. In this case,
$p^r(p^{2^r}-1)\equiv 0 \pmod 5$.

\item[{\rm{(8)}}] A group isomorphic to $\PSL(2,p^m)$ where $m$
divides $r$.

\item[{\rm{(9)}}]  A group isomorphic to $\PGL(2,p^m)$ where $2m$
divides $r$.

\end{itemize}
\end{theorem}
{}From Dickson's classification, all subgroups of $\PGL(2,n)$ with
$n$ odd, can also be obtained, see \cite{maddenevalentini1982}.

Let $n\geq 5$ odd. Then the subgroups listed in {\rm{(1)}} and
{\rm{(2)}} form a partition of $\PSL(2,n)$, that is, every
non-trivial element of $\PSL(2,n)$ belongs exactly one of those
subgroups, see \cite{suzuki1961}. This has the following
corollary.
\begin{proposition}
\label{prop0} Let $n\geq 5$ odd. Any two maximal cyclic subgroups
of $\PSL(2,n)$ have trivial intersection.
\end{proposition}
If $n\geq 5$ is odd, the number of involutions in $\PGL(2,n)$ is
equal to $n^2$.
\begin{proposition}
\label{prop1} Let $n\geq 5$ be odd.
\begin{itemize}
    \item[\rm(I)] $\varphi_{(a,b,c,d)}\in \PGL(2,n)$ is an involution if
and only if $a+d=0$.
    \item[\rm(II)]If $n\equiv 1 \pmod 4$, then $\PSL(2,n)$ has $\ha
n(n+1)$ involutions. Each has exactly exactly two fixed points on
$\Omega$, while no involution in $\PGL(2,n)\setminus \PSL(2,n)$
has a fixed point on $\Omega$.
    \item[\rm(III)] If $n\equiv 3 \pmod 4$, then
$\PSL(2,n)$ has $\ha n(n-1)$ involutions. Each has no fixed point
on $\Omega$, while each involution in $\PGL(2,n)\setminus
\PSL(2,n)$ has exactly two fixed points on $\Omega$.
\end{itemize}
\end{proposition}
\begin{proof} A direct computation shows that $\varphi_{(a,b,c,d)}\in \PGL(2,n)$
is an involution if and only if $b(a+d)=0$ and $c(a+d)=0$. The
latter condition is satisfied  when either  $a+d=0$ or $b=c=0$.
Furthermore, since $\varphi_{(a,0,0,d)}$ is an involution if and
only if $a^2=d^2$ but $a\neq d$, assertion (I) follows.

To show (II) and  (III) take an involution
$\varphi_{(a,b,c,-a)}\in \PGL(2,n)$. A direct computation shows
that $\varphi_{(a,b,c,-a)}$ has two or zero fixed points on
$\Omega$ according as\, $-(a^2-bc)$ is in $\square$ or not. Since
$-1\in \square$ if and only if $n\equiv 1 \pmod 4$, assertions
(II) and (III) follow.
\end{proof}
\begin{proposition}
\label{prop2} Let $n\geq 5$ odd.
\begin{itemize}
    \item[\rm(x)] The elements of $\PGL(2,n)$ of order $p$ are
contained in $\PSL(2,n)$.
    \item[\rm(xx)] Any two elements of $\PSL(2,n)$ of order $p$ are conjugate in
$\PGL(2,n)$.
    \item[\rm(xxx)] The elements of $\PSL(2,n)$ of order $p$ form two different conjugacy classes in $\PSL(2,n)$.
\end{itemize}
\end{proposition}
\begin{proof} In the natural $2$-transitive permutation
representation, the elements $\varphi_{(a,b,c,d)}$ with
$a=d=1,c=0$ and $b\in \Fn$ form a Sylow $p$-subgroup $S_p$ of
$\PGL(2,n)$. Actually, all such elements $\varphi_{(a,b,c,d)}$ are
in $\PSL(2,n)$.

To show (x), it is enough to observe that $\PSL(2,n)$ is
self-conjugate in $\PGL(2,n)$ and that any two Sylow $p$-subgroups
are conjugate in $\PGL(2,n)$.

Take two non-trivial elements in $S_p$, say
$\varphi_1=\varphi_{(1,b,0,1)}$ and
$\varphi_2=\varphi_{(1,b',0,1)}$. Let $a=b'/b$, and
$\varphi=\varphi_{(a,0,0,1)}$. Then
$\varphi_2=\varphi\,\varphi_1\varphi^{-1}$ showing that
$\varphi_2$ is conjugate to $\varphi_1$ in $\PGL(2,n)$. This
proves (xx). Note that if $a\in \square$, then $\varphi_2$ is
conjugate to $\varphi_1$ in $\PSL(2,n)$.

Take any two distinct elements of $\PSL(2,n)$ of order $n$. Every
element of $\PGL(2,n)$ of order $p$ has exactly one fixed point in
$\Omega$ and $\PSL(2,n)$ is transitive on $\Omega$. Therefore, to
show (xxx), we may assume that both elements are in $S_p$. So,
they are $\varphi_1$ and $\varphi_2$ with $b,b'\in \Fn\setminus
\{0\}.$ Assume that $\varphi_2$ is conjugate to $\varphi_1$ under
an element $\varphi\in \PSL(2,n)$. Since $\varphi$ fixes $\infty$,
we have that $\varphi=\varphi_{(a,u,0,1)}$ with $a,u \in F_n$ and
$a\neq 0$. But then $a=b/b'$. Therefore, $\varphi_2$ is conjugate
to $\varphi_1$ under $\PSL(2,n)$ if and only if $b'/b\in \square$.
This shows that $\varphi_1$ and $\varphi_2$ are in the same
conjugacy class if and only $b$ and $b'$ have the same quadratic
character in $\Fn$. This completes the proof.
\end{proof}

\section{The projective unitary group}

Let $\cU$ be the classical unital in $\PG(2,n^2)$, that is, the
set of all self-conjugate points of a non-degenerate unitary
polarity $\Pi$ of $\PG(2,n^2)$. Then $|\cU|=n^3+1$, and at each
point $P\in \cU$, there is exactly one 1-secant, that is, a  line
$\ell_P$ in $\PG(2,n^2)$ such that $|\ell_P\cap \cU|=1$. The pair
$(P,\ell_P)$ is a pole-polar pair of $\Pi$, and hence $\ell_P$ is
an absolute line of $\Pi$. Each other line in $\PG(2,n^2)$ is a
non-absolute line of $\Pi$ and it is an (n+1)-secant of $\cU$,
that is, a line $\ell$ such that $|\ell\cap \cU|=n+1$, see
\cite[Chapter II.8]{hughes-piper1973}.

An explicit representation of $\cU$ in $\PG(2,n^2)$ is as
follows. Let $$ M = \{m\in \fns\mid m^n + m = 0\}.$$ Take an
element $c\in \fns$ such that $c^n+c+1=0$. A homogeneous
coordinate system in $\PG(2,n^2)$ can be chosen so that
$$
\cU = \{X_{\infty}\}\cup \{U=(1,u,u^{n+1}+c^{-1}m)\mid
u\in\fns,\,m\in M\}.
$$
Note that $\cU$ consists of all $\fns$-rational points of the
Hermitian curve of homogeneous equation
$cX_0^nX_2+c^nX_0X_2^n+X_1^{n+1}=0.$

The {\em projective unitary group $\PGU(3,n)$} consists of all
projectivities of $\PG(2,n^2)$ which commute with $\Pi$.
$\PGU(3,n)$ preserves $\cU$ and can be viewed as a permutation
group on $\cU$, since the only projectivity in $\PGU(3,n)$ fixing
every point in $\cU$ is the identity. The group $\PGU(3,n)$ is a
$2$-transitive permutation group on $\Omega$, and this is defined
to be the {\em natural $2$-transitive permutation representation
of $\PGU(3,n).$} Furthermore, $|\PGU(3,n)|=(n^3+1)n^3(n^2-1)$.

With $\mu=\gcd (3,n+1)$, the group $\PGU(3,n)$ contains a normal
subgroup $\PSU(3,n)$, the {\em special unitary group}, of index
$\mu$ which is still a $2$-transitive permutation group on
$\Omega$. This is defined to be the {\em natural $2$-transitive
permutation representation of $\PSU(3,n)$.}

For $n>2$, $\PSU(3,n)$ is a non-abelian simple group, but
$\PSU(3,2)$ is a solvable group.

The maximal subgroups of $\PSU(3,n)$ were determined by Mitchell
\cite{mitchell} for $n$ odd and by Hartley \cite{hartley} for $n$
even, see \cite{hoffer1972}.
\begin{theorem}
\label{hofferclassification} The following is the list of maximal
subgroups of $\PSU(3,n)$ with $n\geq 3$ up to conjugacy$:$
\begin{enumerate}
    \item[\rm(i)] the one-point stabiliser of order $n^3(n^2-1)/\mu;$
    \item[\rm(ii)] the non-absolute line stabiliser of order
$n(n^2-1)(n+1)/\mu;$
    \item[\rm(iii)] the self-conjugate triangle
stabiliser of order $6(n+1)^2/\mu;$
    \item[\rm(iv)] the
normaliser of a cyclic Singer group of order $3(n^2-n+1)/\mu;$
\end{enumerate}
further$,$ for $n=p^k$ with $p>2,$
\begin{enumerate}
    \item[\rm(v)] $\PGL(2,n)$ preserving a conic$;$
    \item[\rm(vi)] $\PSU(3,p^m),$ with $m\mid k$ and $k/m$ odd$;$
    \item[\rm(vii)] the subgroup containing $\PSU(3,p^m)$ as a normal
subgroup of index $3$ when $m\mid k,\ k/m$ is odd$,$ and $3$
divides both $k/m$ and $q+1;$
    \item[\rm(viii)] the Hessian groups of order $216$ when $9\mid(q+1)$, and of order $72$ and $36$ when
$3\mid(q+1);$
    \item[\rm(ix)] $\PSL(2,7)$ when either $p=7$ or
$\sqrt{-7}\not \in \Fq;$
    \item[\rm(x)] the alternating group
$\bA_6$ when either $p=3$ and $k$ is even$,$ or $\sqrt{5}\in \Fq$
but $\Fq$ contains no cube root of unity$;$
    \item[\rm(xi)] the
symmetric group $\bS_6$ for $p=5$ and $k$ odd$;$
    \item[\rm(xii)]
the alternating group $\bA_7$ for $p=5$ and $k$ odd$;$
\end{enumerate}
for $n=2^k,$
\begin{enumerate}
\item[\rm(xiii)] $\PSU(3,2^m)$ with $k/m$ an odd prime$;$
\item[\rm(xiv)] the subgroups containing $\PSU(3,2^m)$ as a normal
subgroup of index $3$ when $k=3m$ with $m$ odd$;$ \item[\rm(xv)] a
group of order $36$ when $k=1$.
\end{enumerate}
\end{theorem}

\begin{proposition}
\label{backgruni} Let $n\geq 3$ be odd. Let $U$ be a cyclic
subgroup of $\PSU(3,n)$ which contains no non-trivial element
fixing a point on $\Omega$. Then $|U|$ divides either $\ha(n+1)$
or $(n^2-n+1)/\mu$.
\end{proposition}
\begin{proof}
Fix a projective frame in $\PG(2,n^2)$ and define the homogeneous
point coordinates $(x,y,z)$ in the usual way. Take a generator $u$
of $U$ and look at the action of $u$ in the projective plane
$\PG(2,\K)$ over the algebraic closure $\K$ of $\fns$. In our
case, $u$ fixes no line point-wise. In fact, if a collineation 
point-wise fixed a line $\ell$ in $\PG(2,\K)$, then $\ell$ would
be a line $\PG(2,n^2)$. But every line in $\PG(2,n^2)$ has a
non-trivial intersection with $\Omega$, contradicting the
hypothesis on the action of $U$.

If $u$ has exactly one fixed point $P$, then $P\in \PG(2,n^2)$ but
$P\not\in \Omega$. Then the polar line $\ell$ of $P$ under the
non-degenerate unitary polarity $\Pi$ is a $(n+1)$-secant of
$\Omega$. Since $\Omega\cap \ell$ is left invariant by $U$, it
follows that $|U|$ divides $n+1$. Since every involution in
$\PSU(3,n)$ has a fixed point on $\Omega$, the assertion follows.

If $u$ has exactly two fixed points $P,Q$, then either $P,Q\in
\PG(2,n^2)$, or $P,Q\in \PG(2,n^4)\setminus \PG(2,n^2)$ and
$Q=\Phi^{(2)}(P),\,P=\Phi^{(2)}(Q)$ where $$\Phi^{(2)}:\,
(x,y,z)\to (x^{n^2},y^{n^2},z^{n^2})$$ is the Frobenius
collineation of $\PG(2,n^4)$ over $\PG(2,n^2)$. In both cases, the
line $\ell$ through $P$ and $Q$ is a line $\ell$ of $\PG(2,n^2)$.
As $u$ has no fixed point in $\Omega$, $\ell$ is not a $1$-secant
of $\Omega$, and hence it is a $(n+1)$-secant of $\Omega$. Arguing
as before shows that $|U|$ divides $\ha(n+1)$.

If $U$ has exactly three points $P,Q,R$, then $P,Q,R$ are the
vertices of a triangle. Two cases can occur according as $P,Q,R\in
PG(2,n^2)$ or $P,Q,R\in \PG(2,n^6)\setminus \PG(2,n^2)$ and
$Q=\Phi^{(3)}(P),\, R=\Phi^{(3)}(Q),\, P=\Phi^{(3)}(R)$ where
$$\Phi^{(3)}:\, (x,y,z)\to
(x^{n^2},y^{n^2},z^{n^2})$$ is the Frobenius collineation of
$\PG(2,n^6)$ over $\PG(2,n^2)$.

In the former case, the line through $P,Q$ is a $(n+1)$-secant of
$\Omega$. Again, this implies that $|U|$ divides $\ha(n+1)$.

In the latter case, consider the subgroup $\Gamma$ of
$\PGL(3,n^2)$, the full projective group of $\PG(2,n^2)$, that
fixes $P,Q$ and $R$. Such a group $\Gamma$ is a Singer group of
$\PG(2,n^2)$ which is a cyclic group of order $n^4+n^2+1$ acting
regularly on the set of points of $\PG(2,n^2)$. Therefore, $U$ is
a subgroup of $\Gamma$. On the other hand, the intersection of
$\Gamma$ and $\PSU(3,n)$ has order $(n^2-n+1)/\mu$, see case (iv)
in Proposition \ref{prop2}.
\end{proof}
\section{The Suzuki group}
A general theory on the Suzuki group is given in \cite[Chapter
XI.3]{huppert-blackburn1982III}.


An {\em ovoid} $\cO$ in $\PG(3,n)$ is a point set with the same
combinatorial properties as an elliptic quadric in $\PG(3,n)$;
namely, $\Omega$ consists of $n^2+1$ points, no three collinear,
such that the lines through any point $P\in\Omega$ meeting
$\Omega$ only in $P$ are coplanar.

In this section, $n=2n_0^2$ with $n_0=n^s$ and $s\geq 1$. Note
that $x^{\varphi}=x^{2q_0}$ is an automorphism of $\Fn$, and
$x^{\varphi^2}=x^2$.

Let $\Omega$ be the {\em Suzuki--Tits ovoid} in $\PG(3,n)$, which
is the only known ovoid in $\PG(3,n)$ other than an elliptic
quadric. In a suitable  homogeneous coordinate system of
$\PG(3,q)$ with $Z_{\infty}=(0,0,0,1)$,
$$
\Omega = \{Z_{\infty}\}\cup\{(1,u,v,uv+u^{2\vf+2}v^{\vf})\mid
u,v\in\Fn\}.
$$

The {\em Suzuki group} $\Sz(n)$, also written $^2B_2(q)$, is the
projective group of $\PG(3,n)$ preserving $\Omega$. The group
$\Sz(n)$ can be viewed as a permutation group on $\Omega$ as the
identity is the only projective transformation in $\Sz(n)$ fixing
every point in $\Omega$. The group $\Sz(n)$ is a $2$-transitive
permutation group on $\Omega$, and this is defined to be the {\em
natural $2$-transitive permutation representation of} $\Sz(n)$.
Furthermore, $\Sz(n)$ is a simple group of order
$(n^2+1)n^2(n-1)$.

The maximal subgroups of $\Sz(n)$ were determined by Suzuki, see
also \cite[Chapter XI.3]{huppert-blackburn1982III}.

\begin{proposition}
The following is the list of maximal subgroups of $\Sz(n)$ up to
conjugacy$:$
    \begin{enumerate} \item[\rm(i)] the one-point
stabiliser of order $n^2(n-1);$ \item[\rm(ii)] the normaliser of a
cyclic Singer group of order $4(n+2n_0+1);$ \item[\rm(iii)] the
normaliser of a cyclic Singer group of order $4(n-2n_0+1);$
\item[\rm(iv)] $\Sz(n')$ for every $n'$ such that $n=n^m$ with $m$
prime.
\end{enumerate}
\end{proposition}
\begin{proposition}
\label{szupart}
 The subgroups listed below form a partition of
$\Sz(n):$
\begin{enumerate}
    \item[\rm(v)] all subgroups of order $n^2;$
    \item[\rm(vi)] all cyclic subgroups
of order $n-1;$
    \item[\rm(vii)] all cyclic Singer subgroups of
order $n+2n_0+1;$
    \item[\rm(viii)] all cyclic Singer subgroups of
order $n-2n_0+1$.
\end{enumerate}
\end{proposition}

\begin{proposition}
\label{backgrsuz} Let $U$ be a cyclic subgroup of $\Sz(n)$ which
contains no non-trivial element fixing a point on $\Omega$. Then
$|U|$ divides either $n-2n_0+1$ or $(n+2n_0+1)$.
\end{proposition}
\begin{proof}
Take a generator $u$ of $U$. Then $u$, and hence $U$, is contained
in one of the subgroups listed in Proposition \ref{szupart}. More
precisely, since $u$ fixes no point, such a subgroup must be of
type (v) or (vi).
\end{proof}
\section{The Ree group}
The Ree group can be introduced in a similar way using the
combinatorial concept of an ovoid, this time in the context of
polar geometries, see for instance \cite[Chapter
XI.13]{huppert-blackburn1982III}.

An {\em ovoid} in the polar space associated to the non-degenerate
quadric $\cQ$ in the space $\PG(6,n)$ is a point set of size
$n^3+1$, with no two of the points conjugate with respect to the
orthogonal polarity arising from $\cQ$.

In this section, $n=3n_0^2$ and $n_0=3^s$ with $s\geq 0$. Then
$x^{\varphi}=x^{3n_0}$ is an automorphism of $\Fn$, and
$x^{\varphi^2}=x^3$.

Let $\Omega$ be the {\em Ree--Tits ovoid} of $\cQ$. In a suitable
homogenous coordinate system of $\PG(6,n)$ with
$Z_{\infty}=(0,0,0,0,0,0,1)$, the quadric is defined by its
homogenous equation $X_3^2+X_0X_6+X_1X_5+X_2X_4=0$, and
$$ \Omega = \{Z_{\infty}\}\cup\{(1,u_1,u_2,u_3,v_1,v_2,v_3)\},$$
with
\begin{eqnarray*}
&& v_1(u_1,u_2,u_3)   =  u_1^2u_2-u_1u_3+u_2^{\vf}-u_1^{\vf+3}, \\
&& v_2(u_1,u_2,u_3)    =
u_1^{\vf}u_2^{\vf}-u_3^{\vf}+u_1u_2^2+u_2u_3-u_1^{2\vf+3},\\
&& v_3(u_1,u_2,u_3) =\\
  &   & \qquad u_1u_3^{\vf}-u_1^{\vf+1}u_2^{\vf}+
u_1^{\vf+3}u_2+u_1^2u_2^2-u_2^{\vf+1}-u_3^2+u_1^{2\vf+4},
 \end{eqnarray*}
for $u_1,u_2,u_3\in\Fn$.

The {\em Ree group} $\Ree(n)$, also written $^2G_2(n)$, is the
projective group of $\PG(6,n)$ preserving $\Omega$. The group
$\Ree(n)$ can be viewed as a permutation group on $\Omega$ as the
identity is the only projective transformation in $\Ree(n)$ fixing
every point in $\Omega$. The group $\Ree(n)$ is a $2$-transitive
permutation group on $\gO$, and this is defined to be the {\em
natural $2$-transitive permutation representation of} $\Ree(n)$.
Furthermore, $|\Ree(n)|=(n^3+1)n^3(n-1)$. For $n_0>1$, the group
$\Ree(n)$ is simple, but $\Ree(3)\cong {\rm{P \Gamma L}}(2,8)$ is
a non-solvable group with a normal subgroup of index $3$.

For every prime $d>3$, the Sylow $d$-subgroups of $\Ree(n)$ are
cyclic, see \cite[Theorem 13.2\,(g)]{huppert-blackburn1982III}.
Put
\begin{eqnarray*}
&& w_1(u_1,u_2,u_3) =  -u_1^{\vf+2}+u_1 u_2- u_3,\\
&&w_2(u_1,u_2,u_3)  =  u_1^{\vf+1}u_2+u_1^{\vf}u_3-u_2^2,\\
&&w_3(u_1,u_2,u_3) =\\
 & & \qquad    u_3^{\vf}+(u_1u_2)^{\vf}-u_1^{\vf+2}u_2-u_1u_2^2+u_2u_3-u_1^{\vf+1}u_3-u_1^{2\vf+3},\\
&&w_4(u_1,u_2,u_3)  =  u_1^{\vf+3}-u_1^2u_2-u_2^{\vf}-u_1u_3.
\end{eqnarray*}
Then a Sylow $3$-subgroup $S_3$ of $\Ree(n)$ consists of the
projectivities represented by the matrices,
$$
\left[\begin{array}{ccccccc} 1 & 0 & 0 & 0 & 0 & 0 & 0 \\ a & 1 &
0 & 0 & 0 & 0 & 0
\\ b & a^{\vf} & 1 & 0 & 0 & 0 & 0
\\ c & b-a^{\vf+1} & -a & 1 & 0 & 0 & 0
\\ v_1(a,b,c) & w_1(a,b,c) & -a^2 & -a & 1 & 0 & 0
\\ v_2(a,b,c) & w_2(a,b,c) & ab+c & b & -a^{\vf} & 1 & 0
\\ v_3(a,b,c) & w_3(a,b,c) & w_4(a,b,c) & c & -b+a^{\vf+1}& -a & 1
\end{array}
\right]
$$
for $a,b,c\in\Fn$. Here, $S_3$ is a normal subgroup of
$\Ree(n)_{Z_{\infty}}$ of order $n^3$ and regular on the remaining
$n^3$ points of $\Omega$. The stabiliser $\Ree(n)_{Z_{\infty},O}$
with $O=(1,0,0,0,0,0,0)$ is the cyclic group of order $n-1$
consisting of the projectivities represented by the diagonal
matrices,
$$
\diag(1,d,d^{\vf+1},d^{\vf+2},d^{\vf+3},d^{2\vf+3}, d^{2\vf+4})$$
for $d\in\Fn$. So the stabiliser $\Ree(n)_{Z_{\infty}}$ has order
$n^3(n-1)$.

The group $\Ree(n)$ is generated by $S_3$ and
$\Ree(n)_{Z_{\infty},O}$, together with the projectivity $W$ of
order $2$ associated to the matrix,
$$
\left[\begin{array}{ccccccc} 0 & 0 & 0 & 0 & 0 & 0 & 1 \\ 0 & 0 &
0 & 0 & 0 & 1 & 0
\\ 0 & 0 & 0 & 0 & 1 & 0 & 0
\\ 0 & 0 & 0 & 1 & 0 & 0 & 0
\\ 0 & 0 & 1 & 0 & 0 & 0 & 0
\\ 0 & 1 & 0 & 0 & 0 & 0 & 0
\\ 1 & 0 & 0 & 0 & 0 & 0 & 0
\end{array}
\right],
$$ that interchanges  $Z_\infty$ and $O$.
Here, $W$ is an involution and it fixes exactly $n+1$ points of
$\Omega$. Furthermore, $\Ree(n)$ has a unique conjugacy classes of
involutions, and hence every involution in $\Ree(n)$ has $n+1$
fixed points in $\Omega$.

Assume that $n=n'^t$ with an odd integer $t=2v+1$, $v\geq 1$. Then
$\Fn$ has a subfield ${\mathbb{F}}_{n'}$, and $\PG(6,n)$ may be
viewed as an extension of $\PG(6,n')$. Doing so, $\cQ$ still
defines a quadric in $\PG(6,n')$, and the points of $\Omega$
contained in $\PG(6,n')$ form an ovoid, the Ree-Tits ovoid of
$\cQ$ in $\PG(6,n')$. The associated Ree group $\Ree(n')$ is the
subgroup of $\Ree(n)$ where the above elements $a,b,c,d$ range
over ${\mathbb{F}}_{n'}$.

The maximal subgroups of $\Ree(n)$ were determined by Migliore
and, independently, by Kleidman \cite[Theorem C]{kleidman1988},
see also \cite[Lemma 3.3]{flp}.
\begin{proposition}
\label{reegroupclassification} The following is the list of
maximal subgroups of $\Ree(n)$ with $n>3$ up to conjugacy$:$
\begin{enumerate}
\item[\rm(i)] the one-point stabiliser of order $n^3(n-1);$

\item[\rm(ii)] the centraliser of an involution $z\in \Ree(n)$
isomorphic to $\langle z \rangle \times \PSL(2,n)$ of order
$n(n-1)(n+1);$

\item[\rm(iii)] a subgroup of order $6(n+3n_0+1),$ the normaliser
of a cyclic Singer group of order $n+3n_0+1;$

\item[\rm(iv)] a subgroup of order $6(n-3n_0+1),$ the normaliser
of a cyclic Singer order of order $6(n-3n_0+1);$

\item[\rm(v)] a subgroup of order $6(n+1),$ the normaliser of a
cyclic subgroup of order $n+1;$

\item[\rm(vi)] $\Ree(n')$ with $n=n'^t$ and $t$ prime.
\end{enumerate}
\end{proposition}

\begin{proposition}
\label{backgrree} Let $U$ be a cyclic subgroup of $\Ree(n)$ with
$n>3$ which contains no non-trivial element fixing a point on
$\Omega$. Then $|U|$ divides either $\ha(n+1)$, or $n-3n_0+1$ or
$n+3n_0+1$.
\end{proposition}
\begin{proof} Every involution in $\Ree(n)$
has exactly $n+1$ fixed points on $\Omega$, and every element in
$\Ree(n)$ whose order is $3$ fixes exactly one point in $\Omega$.
Therefore, neither $3$ nor $2$ divides $|U|$. Furthermore, if $U$
is contained in a subgroup (iii), then $U$ preserves the set of
fixed points of $z$, and hence $|U|$ divides $\ha(n+1)$.

Now, assume that $U$ is contained in a subgroup (iii) or (iv), say
$N$. Let $S$ be the cyclic Singer subgroup of $N$. We show that
$U$ is contained in $S$. Suppose on the contrary that $S\cap U\neq
U$. Then $SU/S$ is a non-trivial subgroup of factor group $N/S$.
Hence either $2$ or $3$ divides $|SU/S|$. Since
$|SU/S|=|S|\cdot|U|/|S\cap U|$ and neither $2$ nor $3$ divides
$|S|$, it follows that either $2$ or $3$ divides $|U|$. But this
is impossible by the preceding result.

If $U$ is contained in a subgroup (v), say $N$, we may use the
preceding argument. Let $S$ be the cyclic subgroup of $N$. Arguing
as before, we can show that $U$ is a subgroup of $S$.

Finally, we deal with the case where $U$ is contained in a
subgroup (vi) which may be assumed to be $\Ree(n')$ with
$$n=n'^{(2v+1)},\,\, v\geq 1;$$ equivalently
$$s=2uv+u+v.$$
 Without loss of generality, $U$ may
be assumed not be contained in any subgroup $\Ree(n'')$ of
$\Ree(n')$.

If $n'=3$ then $U$ is a subgroup of $\Ree(3)\cong {\rm{P\Gamma
L}}(2,8)$. Since $|{\rm{P\Gamma L}}(2,8)|=2^3\cdot 3^3\cdot 7$,
and neither $2$ nor $3$ divides $|U|$, this implies that $|U|=7$.
On the other hand, since $n=3^k$ with $k$ odd, $7$ divides
$n^3+1$. Therefore, $7$ divides $n^3+1=(n+1)(n+3n_0+1)(n-3n_0+1)$
whence the assertion follows.

For $n'>3$, the above discussion can be repeated for $n'$ in place
of $n$, and this gives that $|U|$ divides either $n'+1$ or
$n'+3n_0'+1$ or $n'+3n_0'+1$. So, we have to show that each of
these three numbers must divide either $n+1$, or $n+3n_0+1$, or
$n-3b_0+1$.

If $U$ divides $n'+1$ then it also divide $n+1$ since $n$ is an
odd power of $n'$. For the other two cases, the following result
applies for $n_0=k$ and $n_0'=m$.

\begin{claim}{\rm{\cite[V.~V\'igh]{vigh2008}}}
Fix an $u\geq 0$, and let $m=3^u,\,d^{\pm}=3m^2\pm 3m+1$. For a
non-negative integer $v$, let $s=2uv+u+v,\, k=3^s$, and
$$M_1(v)=3k^2+3k+1,\, M_2(v)=3k^2+1,\, M_3(v)=3k^2-3k+1.$$ Then for
all $v\geq 0$, $d^{\pm}$ divides at least one of $M_1(v)$,
$M_2(v)$ and $M_3(v)$.
\end{claim}
We prove the claim for $d^{+}=d=3m^2+3m+1$, the proof for other
case $d^{-}=m^2-3m+1$ being analog.

We use induction on $v$. We show first that the claim is true for
$v=0,1,2$, then we prove that the claim holds true when stepping
from $v$ to $v+3$.

Since $M_1(0)=d$, the claim trivially holds for $v=0$.

For $v=1$ we have the following equation:
$$(3^{2u+1}+3^{u+1}+1)(3^{4u+2}-3^{3u+2}+3^{2u+2}-3^{2u+1}-3^{u+1}+1)=3^{6u+3}+1,$$
whence
\begin{equation} \label{alap1}
3^{2u+1}+3^{u+1}+1=d\;|\;M_2(1)=3^{6u+3}+1.
\end{equation}
Similarly,
\begin{eqnarray*}(3^{2u+1}+3^{u+1}+1)(3^{8u+4}-3^{7u+4}+3^{6u+4}-3^{6u+3}-3^{5u+3})=\\=3^{10u+5}-3^{5u+3}+1-(3^{6u+3}+1).\end{eqnarray*}
On the other hand, using (\ref{alap1}) we obtain that
\begin{equation*}
3^{2u+1}+3^{u+1}+1=d\;|\;M_3(2)=3^{10u+5}-3^{5u+3}+1,
\end{equation*}
which gives the claim for $v=2$.

Furthermore, using (\ref{alap1}) together with
\begin{eqnarray*}
M_2(v+3)-M_2(v)=(3^{4uv+14u+2v+7}+1)-(3^{4uv+2u+2v+1}+1)=\\
=3^{4uv+2u+2v+1}(3^{6u+3}+1)(3^{6u+3}-1)
\end{eqnarray*}
we obtain that \begin{equation}\label{alap2} d\;|\;
M_2(v+3)-M_2(v).
\end{equation}
Now, direct calculation shows that
$$M_1(v+3)-M_3(v)=M_2(v+3)-M_2(v)+3^{2uv+u+v+1}\cdot M_2(1).$$
{}From (\ref{alap1}) and (\ref{alap2}),
\begin{equation*}
d\;|\; M_1(v+3)-M_3(v).
\end{equation*}
Similarly,
$$M_3(v+3)-M_1(v)=M_2(v+3)-M_2(v)-3^{2uv+u+v+1}\cdot M_2(1),$$
and so
\begin{equation*}
d\;|\; M_3(v+3)-M_1(v).
\end{equation*}
This finishes the proof of the Claim and hence it completes the
proof of Proposition \ref{backgrree}.
\end{proof}
One may ask for a proof that uses the structure of $\Ree(n)$ in
place of the above number theoretic Claim. This can be done as
follows.

Take a prime divisor $d$ of $|U|$. As we have pointed out at the
beginning of the proof of Proposition \ref{backgrree}, $U$ has no
elements of order $2$ or $3$. This implies that $d>3$. In
particular, the Sylow $d$-subgroups of $\Ree(n)$ are cyclic and
hence are pairwise conjugate in $\Ree(n)$.

Since $|U|$ divides $n^3+1$, and $n^3+1$ factorizes into
$(n+1)(n+3n_0+1)(n-3n_0+1)$ with pairwise co-prime factors, $d$
divides just one of this factors, say $v$. From Proposition
\ref{reegroupclassification}, $\Ree(n)$ has a cyclic subgroup $V$
of order $v$. Since $d$ divides $v$, $V$ has a subgroup of order
$d$. Note that $V$ is not contained in $\Ree(n')$ as $v$ does not
divide $|\Ree(n')|$.

Let $D$ be a subgroup of $U$ of order $d$. Then $D$ is conjugate
to a subgroup of $V$ under $\Ree(n)$. We may assume without loss
of generality that $D$ is a subgroup of $V$.

Let $\mathcal{C}(D)$ be the centralizer of $D$ in $\Ree(n)$.
Obviously, $\mathcal{C}(D)$ is a proper subgroup of $\Ree(n)$.
Since both $U$ and $V$ are cyclic groups containing $D$, they are
contained in $\mathcal{C}(D)$. Therefore, the subgroup $W$
generated by $U$ and $V$ is contained in $\mathcal{C}(D)$. To show
that $U$ is a subgroup of $V$, assume on the contrary that the
subgroup $W$ of $\mathcal{C}(D)$ generated by $U$ and $V$ contains
$V$ properly. From Proposition \ref{reegroupclassification}, the
normalizer $\mathcal{N}(V)$ is the only maximal subgroup
containing $V$. Therefore $W$ is a subgroup of $\mathcal{N}(V)$
containing $V$, and $W=UV$. The factor group $W/V$ is a subgroup
of the factor group $\mathcal{N}(V)/V$. From Proposition
\ref{reegroupclassification}, $|W/V|$ divides $6$. On the other
hand,
$$|W/V|=\frac{|U||V|}{|U\cap V||V|}=\frac{|U|}{|U\cap V|}.$$
But then $|U|$ has to divide $6$, a contradiction.

\end{document}